\documentstyle{article}

\input epsf.sty
\input rotate.tex
\def\G{\Gamma}
\def\Om{\Omega}
\def\Si{\Sigma}
\def\T{\Theta}

\def\a{\alpha}
\def\b{\beta}
\def\g{\gamma}
\def\lb{\lambda}
\def\om{\omega}
\def\s{\sigma}
\def\t{\theta}
\def\vp{\varphi}

\def\fl{\forall}
\def\hra{\hookrightarrow}
\def\ify{\infty}

\def\longra{\longrightarrow}
\def\mpo{\mapsto}
\def\op{\oplus}
\def\ot{\otimes}
\def\ov{\overline}
\def\ra{\rightarrow}
\def\sbs{\subset}
\def\ts{\times}
\def\wdg{\wedge}
\def\wt{\widetilde}

\font\tenbb=msbm10
\font\sevenbb=msbm7
\font\fivebb=msbm5
\newfam\bbfam
\textfont\bbfam=\tenbb \scriptfont\bbfam=\sevenbb
\scriptscriptfont\bbfam=\fivebb
\def\bb{\fam\bbfam}

\def\Cb{{\bb C}}
\def\Rb{{\bb R}}
\long\def\comment#1\endcomment{\relax}

\def\Hc{{\cal H}}

\def\Nc{{\cal N}}

\def\Uc{{\cal U}}

\def\build#1_#2^#3{\mathrel{
\mathop{\kern 0pt#1}\limits_{#2}^{#3}}}

\def\diagram#1{\def\normalbaselines{\baselineskip=0pt
\lineskip=10pt\lineskiplimit=1pt}   \matrix{#1}}

\def\hfl#1#2{\smash{\mathop{\hbox to 6mm{\rightarrowfill}}
\limits^{\scriptstyle#1}_{\scriptstyle#2}}}

\def\sevafig#1#2{\centerline{
 \epsfxsize=#2\epsfbox{#1}}}

\font\goth=eufm10

\title{Vanishing of the Kontsevich integrals of the wheels}
\author{Boris Shoikhet}

\begin{document}

\maketitle

\begin{abstract}
We prove that the Kontsevich integrals (in the sense of the 
formality theorem \cite{K}) of all even wheels are equal to 
zero. These integrals appear in the approach to the Duflo 
formula via the formality theorem. The result means that for 
any finite-dimensional Lie algebra {\goth g}, and for 
invariant polynomials $f,g \in [S^{\cdot} (\hbox{\goth 
g})]^{\hbox{\goth g}}$ one has $f \cdot g = f * g$, where 
$*$ is the Kontsevich star-product, corresponding to the 
Kirillov-Poisson structure on $\hbox{\goth g}^*$. We deduce 
this theorem from the result of \cite{FSh} on the 
deformation quantization with traces. 
\end{abstract}

\section{Introduction}\label{sec1}

\subsection{}\label{ssec1.1} First of all, let us recall 
what the Duflo formula is.

Let {\goth g} be a finite-dimensional Lie algebra, we denote 
by $S^{\cdot} (\hbox{\goth g})$ and $\Uc (\hbox{\goth g})$ 
the symmetric and the universal enveloping algebras of the 
Lie algebra {\goth g}, correspondingly. There is the adjoint action 
of {\goth g} on both spaces $S^{\cdot} (\hbox{\goth g})$ and 
$\Uc (\hbox{\goth g})$, and the classical 
Poincar\'e-Birkhoff-Witt map $\vp_{\rm PBW} : S^{\cdot} 
(\hbox{\goth g}) \ra \Uc (\hbox{\goth g})$
\begin{equation}
\vp_{\rm PBW} (g_1 \cdot \ldots \cdot  g_k) = \frac{1}{k!} \ 
\sum_{\s \in \Si_k} \ g_{\s (1)} \ot \ldots \ot g_{\s (k)} 
\label{eq1}
\end{equation}
is an isomorphism of the {\goth g}-modules. In particular, 
it defines a map of invariants
\begin{equation}
[\vp_{\rm PBW}] : [S^{\cdot} (\hbox{\goth g})]^{\hbox{\goth 
g}} \ra [\Uc (\hbox{\goth g})]^{\hbox{\goth g}} \simeq Z 
(\Uc (\hbox{\goth g})) \label{eq2}
\end{equation}
where $Z (\Uc (\hbox{\goth g}))$ is the center of the 
universal enveloping algebra.

The Duflo theorem states that $[S^{\cdot} (\hbox{\goth 
g})]^{\hbox{\goth g}}$ and $Z (\Uc (\hbox{\goth g}))$ are 
isomorphic as algebras, and gives an explicit formula for 
the isomorphism.

For each $k \geq 1$ there exists a canonical invariant 
element ${\rm Tr}_k \in [S^k (\hbox{\goth g})]^*$. It is 
just the trace of $k$-th power of the adjoint action, i.e. 
the symmetrization of the following map ${\rm Tr}_k$:
\begin{equation}
{\rm Tr}_k \, (g) = {\rm Tr}_{\hbox{\goth g}} \, (ad \, (g))^k \, . 
\label{eq3}
\end{equation}
It is easy to show that for semisimple Lie algebras ${\rm Tr}_k = 0$ 
for odd $k$. 
Hence, 
we 
will consider the elements ${\rm Tr}_k$ only for even $k$.

One can consider any element of $[S^k (\hbox{\goth g})]^*$ 
as a differential operator of $k$-th order with constant 
coefficients, acting on $S^{\cdot} (\hbox{\goth g})$. Let us 
note that for a fixed element $\t \in S^{\ell} (\hbox{\goth 
g})$ the values ${\rm Tr}_k (\t)$ are not equal to $0$ only 
for $k \leq \ell$.

Finally, define the map $\vp_{\rm strange} : S^{\cdot} 
(\hbox{\goth g}) \ra S^{\cdot} (\hbox{\goth g})$ by the 
formula
\begin{equation}
\vp_{\rm strange} = \exp \left( \sum_{k \geq 1} \ \a_{2k} \, 
{\rm Tr}_{2k} \right) \label{eq4}
\end{equation}
where the rational numbers $\a_{2k}$ are defined from the 
formula
\begin{equation}
\sum_{k \geq 0} \ \a_{2k} \, q^{2k} = \log \, 
\sqrt{\frac{e^{q/2} - e^{-q/2}}{q}} \, . \label{eq5}
\end{equation}

\medskip

\noindent {\bf Theorem.} (Duflo) {\it For a 
finite-dimensional Lie algebra {\goth g}, the restriction of 
the map $\vp_{\rm PBW} \circ \vp_{\rm strange} : S^{\cdot} 
(\hbox{\goth g}) \ra \Uc (\hbox{\goth g})$ to the space 
$[S^{\cdot} (\hbox{\goth g})]^{\hbox{\goth g}}$ defines an 
isomorphism of the algebras
$$
[\vp_{\rm PBW} \circ \vp_{\rm strange}] : [S^{\cdot} 
(\hbox{\goth g})]^{\hbox{\goth g}} \simeq Z (\Uc 
(\hbox{\goth g})) \, .
$$
}

\medskip

\subsection{}\label{ssec1.2} Here we outline the 
Kontsevich's approach to the Duflo formula via the formality 
theorem (\cite{K}, Sect.~8).

For any Poisson structure on a finite-dimensional vector 
space $V$, i.e. for a bivector field $\a$ on $V$ such that 
$[\a , \a] = 0$, M.~Kontsevich defined a deformation 
quantization of the algebra structure on functions $C^{\ify} 
(V)$.

Any Lie algebra {\goth g} defines the Kirillov-Poisson 
structure on $\hbox{\goth g}^*$. The Poisson bracket of two {\it 
linear} 
functions on $\hbox{\goth g}^*$, i.e. of two elements of 
{\goth g}, is equal to their bracket: $\{ g_1 , g_2 \} := 
[g_1 , g_2]$. This bracket can be extended to $S^{\cdot} 
(\hbox{\goth g})$ by the Leibniz rule. The corresponding 
bivector field in coordinates $\{ x_i \}$ on {\goth g} is
\begin{equation}
\a = \sum_{i,j,k} \ C_{ij}^k \, x_k \ 
\frac{\partial}{\partial x_i} \wdg \frac{\partial}{\partial 
x_j} \label{eq6}
\end{equation}
where $\{ C_{ij}^k \}$ is the structure constants of the Lie 
algebra {\goth g} in the basis $\{ x_i \}$. Finally, the 
bracket of any two functions is $\{ f,g \} = \a \, (df \wdg 
dg)$.

The Kontsevich deformation quantization of this structure 
defines a star-product on $S^{\cdot} (\hbox{\goth g})$, and 
the deformed algebra $(S^{\cdot} (\hbox{\goth g}) , *)$ is 
isomorphic to the universal enveloping algebra $\Uc 
(\hbox{\goth g})$.

\medskip

\noindent {\bf Theorem.} (Kontsevich \cite{K}) {\it There 
exist numbers $W_{2k}$ and $\a'_{2k}$ such that:}
\begin{itemize}
\item[(i)] {\it for any finite-dimensional Lie algebra 
{\goth g} the map
$$
\vp_W = \exp \left( \sum_{k \geq 1} \ W_{2k} \, {\rm 
Tr}_{2k} \right) : S^{\cdot} (\hbox{\goth g}) \ra S^{\cdot} 
(\hbox{\goth g})
$$
defines the isomorphism of the algebras $[\vp_W] : 
[S^{\cdot} (\hbox{\goth g})]^{\hbox{\goth g}} \build 
\longra_{}^{\sim} [S^{\cdot} (\hbox{\goth g}) , 
*]^{\hbox{\goth g}}$ where $*$ is the Kontsevich 
star-product;}
\item[(ii)] {\it the canonical map $\vp_{\a'} : (S^{\cdot} 
(\hbox{\goth g}) , *) \ra \Uc (\hbox{\goth g})$,
$$
\vp_{\a'} (g_1 * \ldots * g_k) = g_1 \ot \ldots \ot g_k
$$
is equal to
$$
\vp_{\a'} =  \vp_{\rm PBW} \circ \exp \left( \sum_{k \geq 1} 
\ \a'_{2k} \, {\rm Tr}_{2k} \right) \, .
$$
}
\end{itemize}

\medskip

As a consequence, we obtain that the coefficients $\a_{2k}$ in 
the Duflo formula are equal to the sum
\begin{equation}
\a_{2k} = W_{2k} + \a'_{2k} \, . \label{eq7}
\end{equation}

The numbers $W_{2k}$ and $\a'_{2k}$ are defined as integrals 
over configuration spaces. They were not computed in \cite{K}. 
The main result is that $W_{2k}$ and $\a'_{2k}$ do not depend 
on the Lie algebra {\goth g}.

The number $W_{2k}$ is the Kontsevich integral corresponding 
to the wheel with $2k$ vertices, see Figure~1.

\bigskip

\vglue 1cm

\sevafig{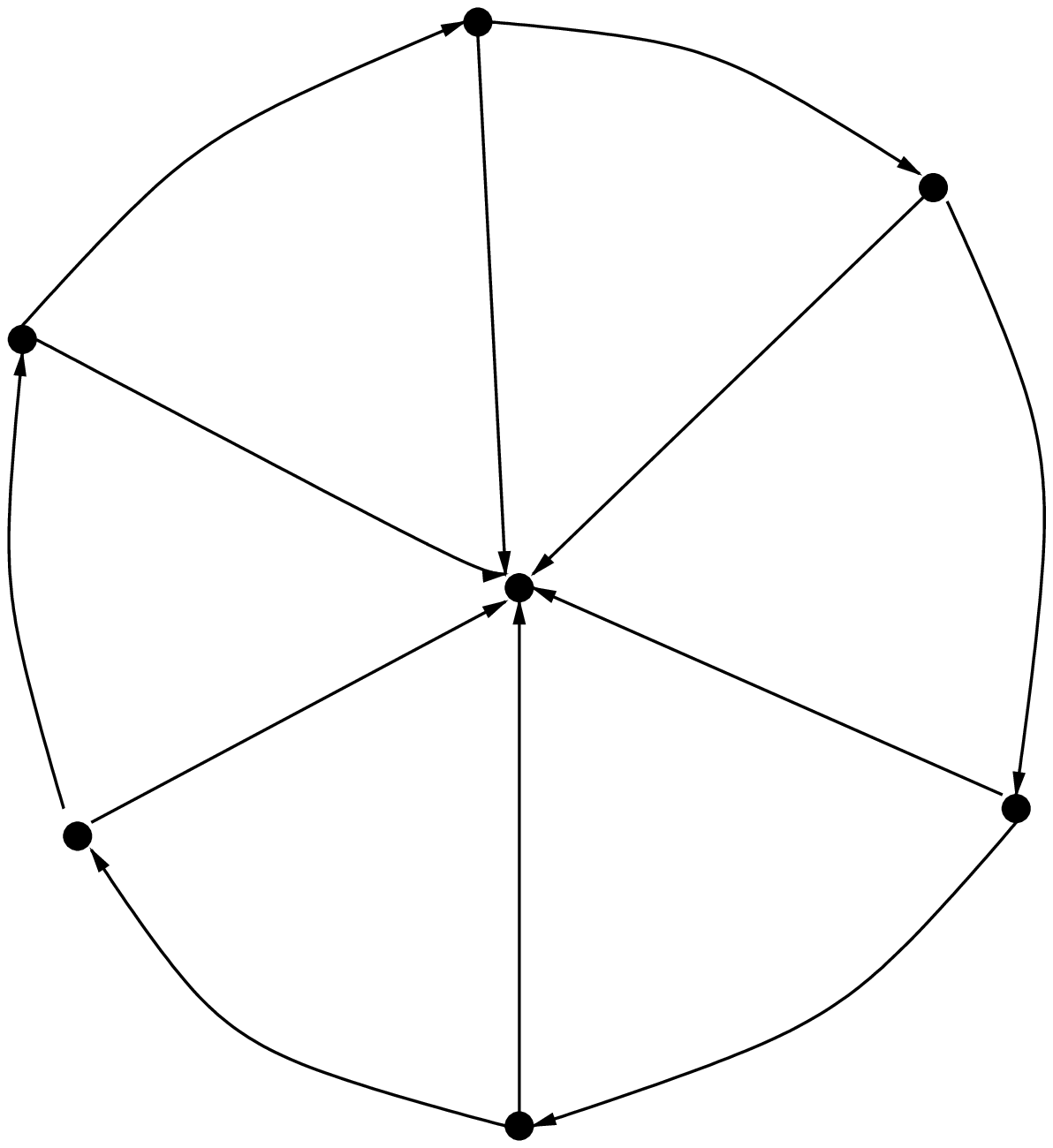}{4cm}

\vglue 1cm

\centerline{}

\smallskip

\centerline{\it The wheel $W_6$}

\bigskip

\noindent {\it Remarks.}

\smallskip

\noindent 1. The Kontsevich star-product on $S^{\cdot} 
(\hbox{\goth g})$ is well-defined only for finite-dimensional 
Lie algebras {\goth g}, while the algebra $\Uc (\hbox{\goth 
g})$ (as well as some other quantizations) is well-defined for 
any (maybe infinite-dimensional) Lie algebra.

\smallskip

\noindent 2. In \cite{K} the Duflo isomorphism was extended 
from the invariants (i.e. zero degree cohomology) to the whole 
 algebras of cohomology. The result is that the map $\vp_{\rm 
PBW} \circ \vp_{\rm strange}$: $S^{\cdot} (\hbox{\goth g}) \ra 
\Uc (\hbox{\goth g})$ defines an isomorphism
$$
[ \vp_{\rm PBW} \circ \vp_{\rm strange} ] : H^{\cdot} 
(\hbox{\goth g} ; S^{\cdot} (\hbox{\goth g})) \ra H^{\cdot} 
(\hbox{\goth g} ; \Uc (\hbox{\goth g})) \, .
$$
In this cohomological setting it seems that an isomorphism 
should exist for infinite-dimensional Lie algebras as well. The 
Duflo formula can not be applied because the traces ${\rm 
Tr}_{2k}$ are ill-defined for infinite-dimensional Lie 
algebras. One of our goals in this work was to understand the 
nature of Duflo isomorphism for infinite-dimensional Lie 
algebras.

\medskip

\subsection{}\label{ssec1.3} The main result of the present 
paper was conjectured by Alberto S. Cattaneo and Giovanni Felder:

\medskip

\noindent {\bf Theorem.} {\it All numbers $W_{2k}$ ($k \geq 
1$) are equal to $0$.
}

\medskip

In an equivalent form, $f \cdot g = f * g$ for any two 
invariant polynomials $f,g \in [S^{\cdot} (\hbox{\goth 
g})]^{\hbox{\goth g}}$ for arbitrary finite-dimensional Lie 
algebra {\goth g} ($*$ is the Kontsevich star-product). Also we 
obtain $\a'_{2k} = \a_{2k}$.

A.~S.~Cattaneo and G.~Felder had computed the number $W_2$ as a
four-dimensional integral and had found that it is equal to 
zero. In this paper we prove that all the numbers $W_{2k}$ are 
equal to $0$ using the deformation quantization with traces 
\cite{FSh}. Let us recall the main result of \cite{FSh}.

Consider a vector space $V$ equipped with a Poisson 
structure $\a$ and a volume form $\Om$ compatible in the 
following way. For any manifold $M$ a volume form on $M$ 
allows to identify polyvector fields with differential forms. 
Then the de~Rham operator on  differential forms defines an 
operator of degree $-1$ on polyvector fields. Such an operator is 
called the 
divergence operator, corresponding to the volume form $\Om$, we 
denote it by ${\rm div}_{\Om}$. It is a second order operator 
with respect to the wedge product of polyvector fields, and 
for any volume form $\Om$, the deffect for the Leibniz rule
\begin{equation}
[ \eta_1 , \eta_2 ] = \pm \, [{\rm div}_{\Om} (\eta_1 \wdg 
\eta_2) - ({\rm div}_{\Om} \, \eta_1) \wdg \eta_2 \pm \eta_1 
\wdg {\rm div}_{\Om} \, \eta_2)] \label{eq8}
\end{equation}
is equal to the Schouten-Nijenhuis bracket of the polyvector 
fields $\eta_1$ and $\eta_2$ and does not depend on $\Om$. As 
a consequence, we obtain
\begin{equation}
{\rm div}_{\Om} \, [ \eta_1 , \eta_2 ] = [ {\rm div}_{\Om} \, 
\eta_1 , \eta_2 ] \pm [\eta_1 , {\rm div}_{\Om} \, \eta_2 ] 
\label{eq9}
\end{equation}
for any volume form $\Om$.

\medskip

\subsubsection{Theorem \cite{FSh}}\label{sssec1.3.1} {\it 
Let $V$ be a finite-dimensional vector space, $\Om$ be a constant volume form
on $V$ (with respect to an affine coordinate system), and
$\a$ be a 
Poisson bivector field on $V$ such that ${\rm div}_{\Om} \, \a = 0$. Then
$$
\int_V f \cdot g \cdot \Om = \int_V (f * g) \cdot \Om
$$
for any two functions $f,g \in C^{\ify} (V)$ one of which has a 
compact support. (Here $*$ is the Kontsevich star-product with 
the harmonic angle function see} \cite{K}, {\it Sect.~$6.2$).}

\medskip

The identity
$$
\int_V f \cdot g \cdot \Om = \int_V (f * g) \cdot \Om
$$
holds for all functions, and we want to remove the integral sign 
when $V = \hbox{\goth g}^*$ and $f,g$ are invariant.

\medskip

\section{Deformation quantization with traces for semi\-simple 
Lie algebras}\label{sec2}

\medskip

\subsection{Theorem.}\label{ssec2.1} {\it Let {\goth g} be a 
semisimple Lie algebra, and let $\a$ be the Poisson-Kirillov 
structure on $\hbox{\goth g}^*$. Then a constant volume form 
$\Om$ on $\hbox{\goth g}^*$ satisfies the equation ${\rm 
div}_{\Om} \, \a = 0$.}

\medskip

\noindent {\it Proof.} For a bivector field
$$
\a = \sum_{i,j} \ a_{ij} \, (x) \, \partial_i \wdg \partial_j 
\qquad (a_{ij} = -a_{ji}) \, ,
$$
its divergence with respect to the constant volume form $\Om = dx_1 
\wdg \ldots \wdg dx_n$ is equal to
\begin{equation}
{\rm div}_{\Om} \, \a = 2 \ \sum_{i,j} \ \partial_i \, (a_{ij} 
\, (x)) \, \partial_j \, . \label{eq10}
\end{equation}
Suppose now that $\a$ is the Poisson-Kirillov bivector field on 
$\hbox{\goth g}^*$. Then, by formula (\ref{eq6}), one has: 
$$
a_{ij} \, (x) = \sum_k \ C_{ij}^k \, x_k
$$
where $C_{ij}^k$ are structure constants of the Lie algebra 
{\goth g} in the basis $\{ x_i \}$. Then, by (\ref{eq10}), the 
condition ${\rm div}_{\Om} \, \a = 0$ is equivalent to
\begin{equation}
\sum_{i,j} \ C_{ij}^i \, \partial_j = 0 \label{eq11}
\end{equation}
or
\begin{equation}
\hbox{for any $j$} \quad \sum_i \ C_{ij}^i = 0 \, . \label{eq12}
\end{equation}
Let us suppose that the Lie algebra {\goth g} is semisimple, and 
the basis $\{ x_i \}$ is chosen in a way compatible with the 
triangular decomposition $\hbox{\goth g} = \Nc_- \op \hbox{\goth h} 
\op 
\Nc_+$. The condition (\ref{eq12}) is nontrivial only when $x_j 
\in \hbox{\goth h}$. Let $e_1 , \ldots , e_{\ell}$ be the positive 
root 
elements in $\Nc_+$, and let $f_1 , \ldots , f_{\ell}$ be the 
dual root elements in $\Nc_-$.

For any element $h \in \hbox{\goth h}$ we set
$$
[h , e_k] = \a_k \, (h) \, e_k
$$
$$
[h , f_k] = \b_k \, (h) \, f_k
$$
where $\a_k , \b_k \in \hbox{\goth h}^*$ are the roots. Then 
(\ref{eq12}) 
holds because $\b_k = -\a_k$.~$\Box$

\bigskip

\noindent {\it Remark.} It follows from the proof that the 
theorem is true also for any unipotent Lie algebra.

\subsection{The vanishing of the wheels.}\label{ssec2.2}

From the previous result and from the result of 
\cite{FSh} (see Sect.~1.3) it follows that for a semisimple Lie 
algebra 
{\goth g} one has
\begin{equation}
\int_{\hbox{\goth g}^*} f \cdot g \cdot \Om =  \int_{\hbox{\goth 
g}^*} (f * g) \cdot \Om \label{eq13}
\end{equation}
where $\Om$ is a constant volume form and $f * g$ is the 
Kontsevich star-product with the harmonic angle function.

We have
\begin{equation}
f * g = f \cdot g + \sum_{k \geq 1} \ \hbar^k \, B_k (f,g) 
\label{eq13bis}
\end{equation}
and
\begin{equation}
B_k (f,g) = (-1)^k \, B_k (g,f) \, . \label{eq14}
\end{equation}
When $f,g$ are invariant, $f * g = g * f$, and, therefore, $B_k 
(f,g) \equiv 0$ for odd $k$. Hence,
\begin{equation}
f * g = f \cdot g + \sum_{k \geq 1} \ \hbar^{2k} \, B_{2k} (f,g) 
\, . \label{eq15}
\end{equation}
Then (13) means that
\begin{equation}
\int_{\hbox{\goth g}^*} B_{2k} (f,g) \cdot \Om =0 
\label{eq16}
\end{equation}
for any $k \geq 1$.

From now on we will work with a real {\it compact} semisimple Lie algebra.
For instance, we replace the complex Lie algebra $sl_n$ to the real Lie algebra $su_n$. We want to apply (16) to invariant $f,g$, that is, to $f$ and $g$ constant on
symplectic leaves. In the compact case $\hbox{\goth g} = Lie G$, where $G$ is a compact Lie group. Therefore, the symplectic leaves, being the orbits of the coajoint action of $G$ on $\hbox{\goth g}^*$, are compact. Hence, we then have some freedom in manipulations with invariant functions with copact support.

We will prove that $W_{2k}$=0 by induction on $k$. Let us suppose that it is proven for $k<l$. Then, by Theorem 1.2(i), we have for invariant $f,g$ :
\begin{equation}
f*g = f \cdot g + {\hbar}^{2l} \cdot (W_{2l} \cdot Tr_{2l} ( f \cdot g) -
W_{2l} \cdot Tr_{2l} (f) \cdot g - W_{2l} \cdot f \cdot Tr_{2l}(g))
+ O({\hbar}^{2l+2})
\label{eq17}
\end{equation}

Formulas (17) and (18) give:
\begin{equation}
W_{2l} \cdot \int_{\hbox{\goth g}^*}(Tr_{2l}(f \cdot g) - Tr_{2l}(f) \cdot g -
f \cdot Tr_{2l}(g)) \cdot \Om = 0
\label{eq18}
\end{equation}

Now we want to prove that the integral in (19) does not vanish for some $f$ and $g$, and, therefore, $W_{2l} = 0$

We have:
\begin{equation}
Tr_{2l}(f) = \sum_{i} a_{i_1...i_{2l}}(f) {\partial}_{i_1} \cdot ... \cdot 
{\partial}_{i_{2l}} (f) 
\label{eq19}
\end{equation}
where $a_{i_1...i_{2l}}$ are constants (depending on the structure constants $C_{ij}^k$). It is clear from (20) that 
\begin{equation}
\int_{\hbox{\goth g}^*}Tr_{2l}(f \cdot g) \cdot \Om = 0 \label{eq20}
\end{equation}
for any $f$ and $g$ with compact support.

Furthemore, it follows from (20) that
\begin{equation}
\int_{\hbox{\goth g}^*}Tr_{2l}(f) \cdot g \cdot \Om =
\int_{\hbox{\goth g}^*}f \cdot Tr_{2l}(g) \cdot \Om \label{eq21}
\end{equation}

Finally, (19) is equivalent to 

\begin{equation}
2 \cdot W_{2l} \cdot \int_{\hbox{\goth g}^*}f \cdot Tr_{2l}(g) \cdot \Om  = 0
\label{eq22}
\end{equation}
which is satisfied for any invariant $f$ and $g$ with compact support.

For the Lie algebra $su_n$, $n \gg 0$, one can choose an invariant $g$ with a compact support such that $Tr_{2l}(g) \ne 0$. It follows from the general description of the algebra of invariants $[S^{\cdot}(\hbox{\goth g})]^{\hbox{\goth g}}$ for $\hbox{\goth g} = su_n$. Then in (23) $f$ is arbitrary, and we obtain $W_{2l}=0$.

\subsection{Remarks.}\label{ssec2.3}
\subsubsection{Remark.}\label{sssec2.3.1}

It is interesting to note that the integrals of the wheels with the opposite direction of the central arrows are {\it not} equal to $0$. Let us denote by $W_k^{\vee}$ such a wheel. 

\bigskip

\vglue 1cm

\sevafig{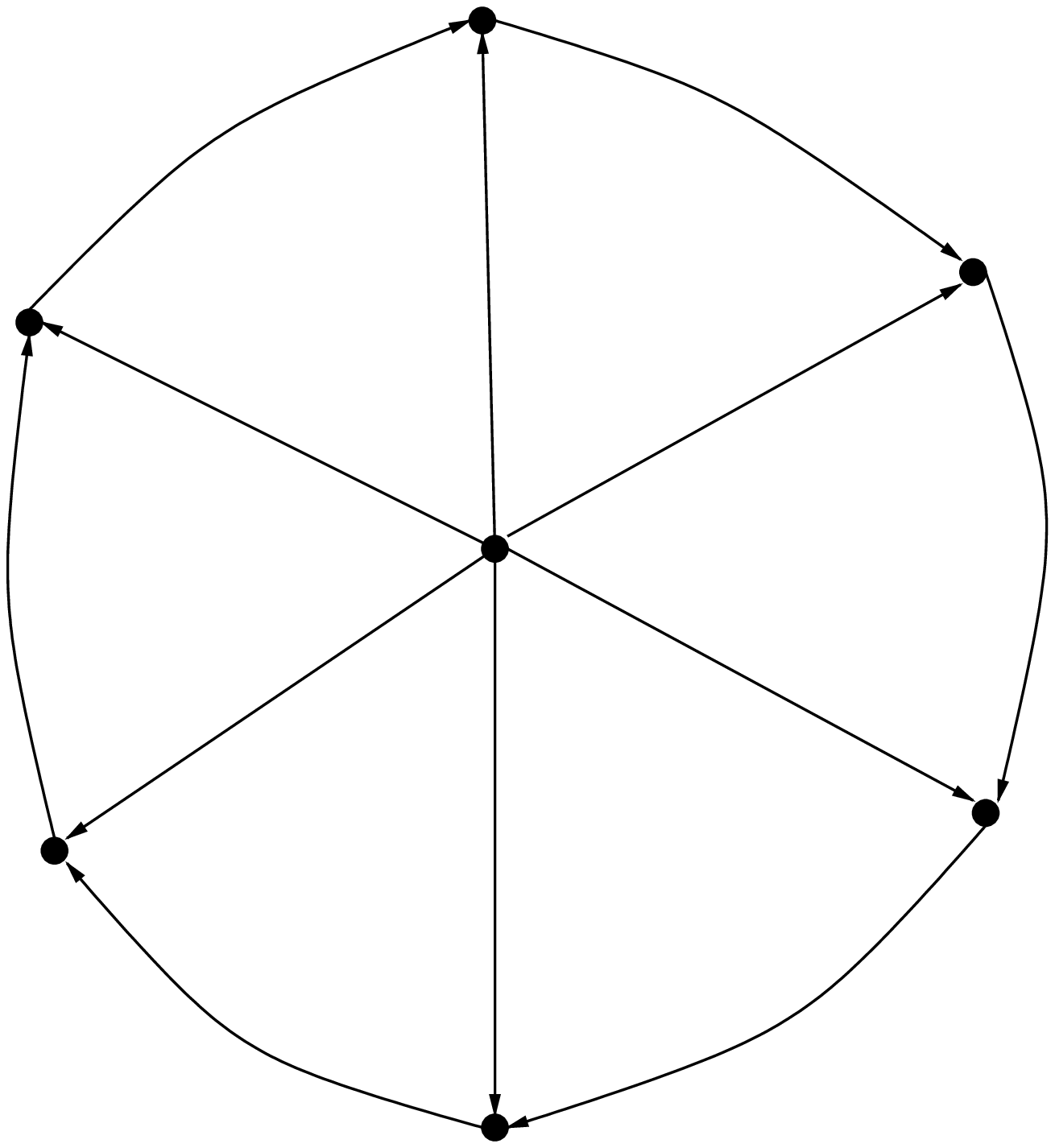}{5cm}

\vglue 1cm

\centerline{}

\smallskip

\centerline{\it The wheel $W_6^{\vee}$}

\bigskip 

Then $W_{2k+1}^{\vee}=0, k \geq 1$, and
$W_{2k}^{\vee}={\alpha}_{2k}$ (see formula (5)), $k \geq 1$.
It follows from an alternative approach to the Duflo formula via the formality theorem, developed in an (unpublished) joint paper of the author with M.Kontsevich.

\subsubsection{Non-linear Poisson structures.}\label{sssec2.3.2}

\noindent {\bf Conjecture} {\it For any Poisson structure $\alpha$  on a vector spase $V$ and for any two functions $f,g \in Z(C^{\infty}(V), *)$ in the center of the deformed algebra one has
$$
f*g=f \cdot g
$$
Here $*$ is the Kontsevich deformation quantization with the harmonic angle function.}

\subsubsection{Tangential deformation quantization.}\label{sssec2.3.3}

Note that our theorem, $f*g=f \cdot g$ for invariant $f,g$, seems to be very closed to the condition of the tangential deformation quantization [CGR], which is $f*g=f \cdot g$ for invariant (=constant on leaves) $f$ and for {\it any} $g$. The last condition means geometrically that we quantize each symplectic leaf separately and then "glue" all the quantized leaves. More algebraically, it means, that in formula (14) all the differential operators $B_k(f,g)$ are formed by a composition of vector fields tangential to the leaves.
 
It turns out, however, that this condition is much stronger; in particular, it is proven in [CGR] that any such a quantization does not exist for semisimple Lie algebras.

\comment
 
We are going to show now that $B_2 (f,g) = 0$ for invariant $f$ 
and $g$. It is clear that a nonzero contributions to $B_2 (f,g)$ 
have the following 4 graphs:

\bigskip

\vglue 1cm

\sevafig{w2.ps}{10cm}

\vglue 1cm

\centerline{Figure 2}

\bigskip

\noindent (see \cite{K} for the definitions). Let us denote the 
graphs showed in Figure 2 by $\G_1$, $\G_2$, $\G_3$, $\G_4$, 
corresp. Let us denote by $w \, (\G)$ the Kontsevich weight of 
a graph $\G$.

\medskip

\noindent {\bf Lemma.}
$$
w \, (\G_1) - w \, (\G_2) - w \, (\G_3) + w \, (\G_4) = 0 \, .
$$

latex
\medskip

\noindent {\it Proof.} It is clear that
$$
\int_{\hbox{\goth g}^*} B_2 (f,g) \cdot \Om = (w \, (\G_1) - w 
\, (\G_2) - w \, (\G_3) + w \, (\G_4)) \ts \int_{\hbox{\goth 
g}^*} f \cdot \T \, (g) \cdot \Om
$$
where $\T$ is the differential operator corresponding to the 
graph $\G$, showed in Figure~3.

\bigskip

\vglue 1cm

\sevafig{w3.ps}{5cm}

\vglue 1cm

\centerline{Figure 3}

\bigskip

On the other hand, it follows from (\ref{eq13}) that 
$\int_{\hbox{\goth g}^*} B_2 (f,g) \cdot \Om = 0$. We have
$$
(w \, (\G_1) - w \, (\G_2) - w \, (\G_3) + w \, (\G_4)) \ts 
\int_{\hbox{\goth g}^*} f \cdot \T \, (g) \cdot \Om = 0
$$
for any $f$. It is clear that $\T \, (g)$ is a non-zero function, so 
Lemma is proven.~$\Box$

\bigskip

Now we are going to show that for invariant $f,g$ one has:
\begin{equation}
B_2 (f,g) = (w \, (\G_1) - w \, (\G_2) - w \, (\G_3) + w \, 
(\G_4)) \ts \Uc_{\G_4} (f,g) \label{eq16}
\end{equation}
where $\Uc_{\G_4}$ is corresponding to the fourth graph on 
Figure~2. For example, let us show that $\Uc_{\G_2} (f,g) = - 
\Uc_{\G_4} (f,g)$ for invariant $f,g$.

\bigskip

\vglue 1cm

\sevafig{w4.ps}{7cm}

\vglue 1cm

\centerline{Figure 4}

\smallskip

\centerline{\it Graph $\G_2$}

\bigskip

One has (see Figure~4):
\begin{equation}
\Uc_{\G_2} (f,g) = \sum_{i,j,k,\ell} \ \partial_{k} \, \a^{ij} 
\cdot \a^{k \ell} \cdot \partial_i (f) \cdot \partial_j \, 
\partial_{\ell} (g) \, . \label{eq17}
\end{equation}
If $g$ is invariant,
\begin{equation}
[x_k , g] = 0 \quad \hbox{for any} \ x_k \label{eq18}
\end{equation}
or
\begin{equation}
\sum_{\ell ,s} \ C_{k\ell}^s \, x^s \, \partial_{\ell} \, g = 0 
\qquad \fl \, k \, . \label{eq19}
\end{equation}
Formula (\ref{eq19}) can be rewritten as
\begin{equation}
\sum_{\ell} \ \a^{k \ell} \, \partial_{\ell} (g) = 0 \qquad 
\hbox{for fixed} \ k \, . \label{eq20}
\end{equation}
Therefore,
\begin{equation}
\sum_{i,j,k,\ell} \ \partial_j \, (\partial_k \, \a^{ij} \cdot 
\a^{k\ell} \cdot \partial_i (f) \cdot \partial_{\ell} (g)) 
\equiv 0 \, . \label{eq21}
\end{equation}
Since $\a^{ij}$ is linear function, by the Leibniz rule, one 
has
\begin{eqnarray}
0 &= &\sum_{i,j,k,\ell} \ \partial_k \, \a^{ij} \cdot \partial_j 
\, \a^{k\ell} \cdot \partial_i (f) \cdot \partial_{\ell} (g) 
\nonumber \\
&+ &\sum_{i,j,k,\ell} \ \partial_k \, \a^{ij} \cdot \a^{k\ell} 
\cdot \partial_i \, \partial_j (f) \cdot \partial_{\ell} (g) + 
\Uc_{\G_2} (f,g) \, . \label{eq22}
\end{eqnarray}
The second summand is equal to $0$ because $\a^{ij}$ is 
skew-symmetric in $i,j$ and $\partial_i \, \partial_j (f)$ is 
symmetric in $i,j$. So, we obtain $\Uc_{\G_2} (f,g) = - 
\Uc_{\G_4} (f,g)$ for invariant $g$.

Analogously, for invariant $f$ one has $\Uc_{\G_3} (f,g) = - 
\Uc_{\G_4} (f,g)$, and for invariant $f,g$ one has $\Uc_{\G_1} 
(f,g) = \Uc_{\G_4} (f,g)$. Therefore, it follows from Lemma 
above that $B_2 (f,g) = 0$ for invariant $f,g$.

Now we want to prove that $W_2 = 0$. It follows from 
Theorem~\ref{ssec1.2}~(i) that if $W_2 \ne 0$ and $B_2 (f,g) = 
0$, then ${\rm Tr}_2$ is a derivation of the algebra $[S^{\cdot} 
(\hbox{\goth g})]^{\hbox{\goth g}}$. It is not true for 
$\hbox{\goth g} = s\ell_2$. So, we are done.

\medskip

\noindent {\it Remark.} It is not clear how to prove that $B_4 (f,g) 
= 0,...$
using 
this method. The cause is that the 
combinatorics of the second part of our proof becomes very 
complicated, and it is not clear how to deduce explicitly the 
identity $B_{2k} (f,g) = 0$ from formula (\ref{eq13}). In the 
next section we develop a different technics based on the orbit 
method.

\section{Vanishing of the wheels}\label{sec3}

\subsection{}\label{ssec3.1} 

Let $\a$ be a Poisson structure on $\Rb^d$ satisfying the following 
two conditions:
\begin{itemize}
\item[(A)] All the symplectic leaves of the structure $\a$ are 
compact;
\item[(B)] If $f,g \in C^{\ify} (\Rb^d)$ are constant along the 
symplectic leaves, $f * g$ also is (here $*$ is the Kontsevich 
deformation quantization with the harmonic angle function).
\end{itemize}

\medskip

\noindent {\it Remarks.} 
\begin{itemize}
\item[(i)] Condition (B) is satisfied in the case of a linear 
structure $\a$, i.e. the Kostant-Kirillov Poisson structure on 
$\hbox{\goth g}^*$, where {\goth g} is a Lie algebra. It follows, for 
example, from Theorem~1.2~(i). The author does not know any other 
examples (except constant and linear cases) when the condition (B) is 
satisfied.
\item[(ii)] For a complex semisimple Lie algebra {\goth g} the 
condition (A) is not satisfied. Fortunately, any such Lie algebra has 
a {\it compact} form $\hbox{\goth g}^{\Rb}$, which is a semisimple 
Lie algebra over $\Rb$ such that $\hbox{\goth g}^{\Rb} \otimes_{\Rb} 
\Cb = \hbox{\goth g}$, and $\hbox{\goth g}^{\Rb} = {\rm Lie} \, (G)$ 
where $G$ is a {\it compact} Lie group. The symplectic leaves are the 
orbits of the coadjoint action of the group $G$. Therefore, they are 
compact, and the condition (A) is satisfied. For our purposes it is 
sufficient to consider the real Lie algebra $\hbox{\goth g}^{\Rb}$ 
instead of the complex Lie algebra {\goth g}.
\end{itemize}

\medskip

\noindent {\bf Theorem.} {\it Let $\a$ be a Poisson structure on 
$\Rb^d$ such that ${\rm div}_{\Om} \, \a = 0$ where $\Om$ is a 
constant volume form on $\Rb^d$. Let us suppose that the conditions 
(A) and (B) above are satisfied. Then $f * g = f \cdot g$ for any two 
functions $f,g \in C^{\ify} (\Rb^d)$ constant along the symplectic 
leaves.}

\medskip

\noindent {\it Proof.} Let $\Uc : \Rb^d \ra \Rb^d$ be the 
transformation, generated by a vector field $v$ on $\Rb^d$, tangent 
to the leaves of the symplectic foliation corresponding to the 
Poisson structure $\a$.

\smallskip

By the Theorem~1.3.1, one has:
\begin{equation}
\int_{\Rb^d} f * g \cdot \Om = \int_{\Rb^d} f \cdot g \cdot \Om \, . 
\label{eq24}
\end{equation}
We apply the formula of the change of variables in integral to the 
transformation $\Uc$. We get:
\begin{equation}
\int_{\Rb^d} (f * g) (\Uc (X)) \cdot \det \Uc (X) \cdot \Om = 
\int_{\Rb^d} f (\Uc (X)) \cdot g (\Uc (X)) \cdot \det \Uc (X) \cdot 
\Om \, . \label{eq25}
\end{equation}
We have: $f (\Uc (X)) \equiv f(x)$ and $g(\Uc (X)) \equiv g(X)$ 
because $f,g$ are constant along the leaves and the transformation 
$\Uc$ preserves the leaves. When the condition (B) is satisfied, we 
also have
$$
(f*g) (\Uc (X)) = (f*g) (x) \, .
$$
In this case, formula (\ref{eq25}) gives:
\begin{equation}
\int_{\Rb^d} (f * g) \cdot \det \, \Uc (X) \cdot \Om = \int_{\Rb^d} f 
\cdot g \cdot \det \, \Uc (X) \cdot \Om \, . \label{eq26}
\end{equation}

This formula is satisfied only when one of the functions $f,g$ has a 
compact support. When the symplectic leaves are not compact, there 
does not exist any constant along the leaves function, except the 
functions with the support on degenerated leaves. When the condition 
(A) is satisfied, any constant along the leaves function can be 
represented as an (infinite) sum of such functions with compact 
support.

\smallskip

The function $\det \, \Uc (X)$ in the domain where $\a$ is regular 
can be (almost) arbitrary. One can suppose that $\Uc \equiv 1$ 
outside a domain in $\Rb^d$, and substract (\ref{eq25}) from 
(\ref{eq26}). We have:
\begin{equation}
\int_{\Rb^d} (f * g) \cdot (\det \, \Uc (X) - 1) \cdot \Om = 
\int_{\Rb^d} f \cdot g \cdot (\det \, \Uc (X) - 1) \cdot \Om \, . 
\label{eq27}
\end{equation}
Here the function $\det \, \Uc (X) - 1$ is (almost) arbitrary 
function with a compact support. We conclude that $f * g = f \cdot 
g$.~$\Box$

\bigskip

\subsection{Vanishing of the wheels}\label{ssec3.2}

Let $k$ be a minimal integral number such that $W_{2k} \ne 0$. 
Then, by Theorem~\ref{ssec1.2}~(i) and by Theorem~\ref{ssec3.4}, 
the element ${\rm Tr}_{2k} \in [S^{2k} (\hbox{\goth 
g}^*)]^{\hbox{\goth g}}$ should be a derivation of the algebra 
$[S^{\cdot} (\hbox{\goth g})]^{\hbox{\goth g}}$ for any semisimple 
Lie algebra {\goth g}. One can show that it is not true for Lie 
algebra $s\ell_n$ for sufficiently large $n$. (The same argument was 
used by Maxim Kontsevich in his proof of Duflo theorem, see 
\cite{K}, Sect.~8.3.4.) Theorem~\ref{ssec1.3} is proven.

\bigskip

\subsection{Remarks.}\label{ssec3.3} 1. The integrals, corresponding
to the wheels, do not depend on a choice of the angle function 
(\cite{K}, Sect.~6.2 and Sect.~8) We have proved that $W_{2k} = 0$ 
for the harmonic angle function (since the result of \cite{FSh} holds 
for harmonic angle function). Hence, it is true for any angle 
function.

\smallskip

\noindent 2. The odd wheels $W_{2k+1}$ are equal to zero by a much 
more simple reasons: it follows from the fact that for the 
harmonic angle function the propagator is invariant under the map 
$p : z \mpo - \ov z$ ($p : \Hc \ra \Hc$, where $\Hc$ is the 
complex upper half-plane), as it was noted in \cite{K}. The same 
argument shows that $B_k (f,g) = (-1)^k \, B_k (g,f)$ for $k \geq 
0$, {\it any} $f,g \in C^{\ify} (V)$ and {\it any} Poisson 
bivector field (not necessarily linear) on a vector space $V$ (see 
\cite{CF}).

\bigskip

\subsection{Deformation quantization with traces and tangential 
deformation quantization}\label{ssec3.4} 

Here we explain why the result of \cite{FSh} may not be true for any 
(compatible with the Poisson structure) volume form, and for the 
Kontsevich star-product.

\smallskip

Let us suppose that we have
\begin{equation}
\int_{\Rb^d} f * g \cdot \wt \Om = \int_{\Rb^d} f \cdot g \cdot \wt 
\Om \label{eq28}
\end{equation}
for any volume form $\wt \Om$ such that ${\rm div}_{\wt \Om} \, \a = 
0$, and for the Kontsevich star-product. Suppose that, in particular, 
${\rm div}_{\Om} \, \a = 0$ for a constant volume form $\Om$ (as it 
is in the semisimple linear case). Then the volume form $\wt \Om = 
\psi \cdot \Om$ satisfies ${\rm div}_{\wt \Om} \, \a = 0$ if and only 
if $\a (d\psi) = 0$, i.e. $\psi$ is constant along the symplectic 
leaves of the Poisson structure $\a$. Then, for any such $\psi$, 
(\ref{eq28}) gives:
\begin{equation}
\int_{\Rb^d} f * g \cdot \psi \cdot \Om = \int_{\Rb^d} f \cdot g \cdot 
\psi \cdot \Om \label{eq29}
\end{equation}
(we suppose that $f$ has a compact support). The left-hand side of 
(\ref{eq29}) is equal to
\begin{equation}
\int_{\Rb^d} f * g \cdot \psi \cdot \Om = \int (g * \psi) \cdot f 
\cdot \Om \label{eq30}
\end{equation}
(apply (\ref{eq28}) for $\wt \Om = \Om$ and the associativity $(f*g) * 
\psi = f * (g*\psi)$). It follows from (\ref{eq29}) and (\ref{eq30}) 
that
\begin{equation}
g * \psi = g \cdot \psi \label{eq31}
\end{equation}
for any $g$ and $\psi$ constant along the symplectic leaves. One can 
deduce 
\begin{equation}\psi * g = \psi \cdot g \label{eq32}
\end{equation}
as well.

\smallskip

A deformation quantization for which (\ref{eq31}) and (\ref{eq32}) are 
satisfied, is called tangential. It is proven in \cite{CGR} that the 
Kontsevich-Kirillov Poisson structure on $\hbox{\goth g}^*$ does not 
admit any {\it tangential} deformation quantization by bidifferential 
operators for any semisimple Lie algebra {\goth g}. Therefore, 
(\ref{eq28}) does not hold.

\subsection{Remark}\label{ssec3.5}

It is very interesting to note that although the Kontsevich integrals 
of the wheels $W_k$ are equal to zero, they do not vanish for very 
similar graphs. Namely, let us denote by $W_k^v$ the wheel obtained 
from the wheel $W_k$ by the reversing of the direction of all central 
edges (they are incoming to the central vertex in the case of $W_k$ 
and they are outcoming in the case of $W_k^v$).

\bigskip

\vglue 1cm

\sevafig{w4bis.ps}{7cm}

\vglue 1cm

\centerline{Figure 5}

\smallskip

\centerline{\it The wheel $W_6^v$}

\bigskip

\noindent It follows from (an unpublished) joint work of the author 
with M.~Kontsevich that the integral of the wheel $W_k^v$ is equal to 
$\a_k$ (see Section~1.1) for even $k$ and to zero for odd $k$.

\section{Duflo formula for infinite-dimensional Lie 
algebras}\label{sec4}

Let {\goth g} be a finite-dimensional Lie algebra. In \cite{K}, 
Sect.~8 M.~Kontsevich proved that the map $\vp_{\rm PBW} \circ 
\vp_{\rm strange} : S^{\cdot} (\hbox{\goth g}) \ra \Uc 
(\hbox{\goth g})$ induces an isomorphism of the {\it algebras} 
$$
[\vp_{\rm PBW} \circ \vp_{\rm strange}] : H^{\cdot} (\hbox{\goth 
g} ; S^{\cdot} (\hbox{\goth g})) \ra H^{\cdot} (\hbox{\goth g} ; 
\Uc (\hbox{\goth g}))
$$
(see Section~\ref{ssec1.2}); the classical Duflo formula is the 
restriction of this theorem to $0$-cohomology part.

It is natural to conjecture analogous statement for arbitrary Lie 
algebras.

\medskip

\noindent {\bf Conjecture.} {\it For any (maybe 
infinite-dimensional) Lie algebra {\goth g} the algebras 
$H^{\cdot} (\hbox{\goth g} ; S^{\cdot} (\hbox{\goth g}))$ and 
$H^{\cdot} (\hbox{\goth g} ; \Uc (\hbox{\goth g}))$ are 
isomorphic.}

\medskip

The map $\vp_{\rm strange}$ does not make sense for an 
infinite-dimensional Lie algebra {\goth g}, because the elements 
${\rm Tr}_{2k} \in [S^{2k} (\hbox{\goth g}^*)]^{\hbox{\goth g}}$ 
are not well-defined, as traces of operators on an 
infinite-dimensional space. Therefore, one needs some 
``regularization'' of these traces. Roughly speaking, one should 
represent each element ${\rm Tr}_{2k}$ as a well-defined quantity 
modulo coboundaries.

It turns out that, in a sense, our result on vanishing 
of the wheels allows to do it.

Let us recall (see Section~\ref{ssec1.2}) that the map $\vp_{\rm 
PBW} \circ \vp_{\rm strange}$ is equal to the composition
\begin{equation}
\diagram{S^{\cdot} (\hbox{\goth g}) &\hfl{\vp_W}{} &(S^{\cdot} 
(\hbox{\goth g}) , *) &\hfl{\vp_{\rm PBW} \circ \vp_{\a'}}{} &\Uc 
(\hbox{\goth g}) \, .} \label{eq33}
\end{equation}
We have proved that $\vp_W = {\rm id}$. On the other hand, the 
algebras $(S^{\cdot} (\hbox{\goth g}) , *)$ and $\Uc (\hbox{\goth 
g})$ are canonically isomorphic. Therefore, the algebras 
$H^{\cdot} (\hbox{\goth g} ; (S^{\cdot} (\hbox{\goth g}) , *))$ 
and $H^{\cdot} (\hbox{\goth g} ; \Uc (\hbox{\goth g}))$ are 
isomorphic as well.

The problem now is that the Kontsevich deformation quantization 
$(S^{\cdot} (\hbox{\goth g}) , *)$ does not exist for an 
infinite-dimensional Lie algebra {\goth g}.

It turns out, however, that the algebra of invariants $[S^{\cdot} 
(\hbox{\goth g}) , *]^{\hbox{\goth g}}$, or, more generally, the 
algebra $H^{\cdot} (\hbox{\goth g} ; (S^{\cdot} (\hbox{\goth g}) , 
*))$ still exists. In other words, in the infinite-dimensional 
case, one has a {\it third} algebra $H_{\rm reg}^{\cdot} 
(\hbox{\goth g} ; (S^{\cdot} (\hbox{\goth g}) , *))$ besides the 
algebras $H^{\cdot} (\hbox{\goth g} ; S^{\cdot} (\hbox{\goth g}))$ 
and $H^{\cdot} (\hbox{\goth g} ; \Uc (\hbox{\goth g}))$.

We have already seen in Section~\ref{ssec2.2} that for  
invariant $f,g \in S^{\cdot} (\hbox{\goth g})$ one has: 
$\Uc_{\G_2} (f,g) = - \Uc_{\G_4} (f,g)$ (see Figure~5)

\bigskip

\vglue 1cm

\sevafig{w5.ps}{9cm}

\vglue 1cm

\centerline{Figure 6}

\smallskip

\centerline{\it An identity for invariant $f,g$}

\bigskip

This identity is the simplest example of the regularization. For 
an infinite-dimensional vector space $V$ the polydifferential 
operator
$$
\Uc_{\G} \, (\g_1 , \ldots , \g_m) \qquad (\g_1 , \ldots , \g_m 
\in T_{\rm poly}^{\cdot} (V))
$$
is well-defined only when the graph $\G$ does not contain any
 oriented cycles between the 
verticies of the first type (see \cite{K}, Sect.~6 for the 
definitions). The procedure described in Section~\ref{ssec2.2} can 
be applied for a regularization of any graph with oriented cycles. 
It means that if a graph $\G$ contains an oriented cycle,  
for invariant $f,g \in S^{\cdot} (\hbox{\goth g})$:
\begin{equation}
\Uc_{\G} (f,g) = \sum_{\wt \G} \ \a_{\wt \G} \, \Uc_{\wt \G} \, 
(f,g) \label{eq34}
\end{equation}
where the graphs $\wt \G$ do not contain any oriented cycles. The 
identity (\ref{eq34}) can be easily generalized for higher 
cohomology. The main problem here is that the regularization is 
not uniquely defined, i.e. the answer ``$\Uc_{\G}^{\rm reg} 
(f,g)$'' depends on the regularization. In the 
finite-dimensional case, however, the answer is independent. 
 This fact motivates the following

\medskip

\noindent {\bf Conjecture.} {\it The regularized algebra $H_{\rm 
reg}^{\cdot} (\hbox{\goth g} ; (S^{\cdot} (\hbox{\goth g}) , *))$ 
does not depend on the regularization.}

\medskip

Finally, the following conjecture is also highly nontrivial:

\medskip

\noindent {\bf Conjecture.} {\it The algebra $H_{\rm reg}^{\cdot} 
(\hbox{\goth g} ; (S^{\cdot} (\hbox{\goth g}) , *))$ is isomorphic 
to the algebra $H^{\cdot} (\hbox{\goth g} ;$ $\Uc (\hbox{\goth g}))$ 
for an infinite-dimensional Lie algebra {\goth g}.}

\medskip

We suppose, by analogy with the finite-dimensional case, 
that the identity map gives an isomorphism ${\rm id} = \vp_W : 
H^{\cdot} (\hbox{\goth g} ; S^{\cdot} (\hbox{\goth g})) \ra 
H^{\cdot} (\hbox{\goth g} ; (S^{\cdot} (\hbox{\goth g}) , *))$. 
Hence, we obtain an isomorphism between the algebras $H^{\cdot} 
(\hbox{\goth g} ; S^{\cdot} (\hbox{\goth g}))$ and $H^{\cdot} 
(\hbox{\goth g} ; \Uc (\hbox{\goth g}))$. It would be very 
interesting to develop this approach on concrete examples.

\bigskip

\endcomment

\bigskip

\noindent {\bf Acknowledgments.} The author is grateful to Boris Feigin and to
Giovanni Felder for inspiring 
discussions. I am thankful to Mme C\'ecile Gourgues for
the quality typing of this text.

\vglue 2cm

\noindent Boris Shoikhet

\noindent Department of Mathematics

\noindent Utrecht University

\noindent P.O. Box 80010

\noindent 3508 TA Utrecht

\noindent The Netherlands

\smallskip

\noindent e-mail address:  borya@mccme.ru

\end{document}

\section{Orbit method and vanishing of the wheels}\label{sec3}

\subsection{}\label{ssec3.1} We start with the following result:

\medskip

\noindent {\bf Key-lemma.} {\it Let $\a$ be a Poisson bivector 
field on a manifold $M$, and let $\Om$ be a volume form on $M$ 
such that ${\rm div}_{\Om} \, \a = 0$. Let $\psi$ be a non-vanishing 
function 
on $M$ which is equal to a constant on each leaf of the 
symplectic foliation on $M$, defined from bivector field $\a$ 
(the constant may depend on the leaf). Then the volume form 
$\wt{\Om} = \psi \cdot \Om$ also satisfies ${\rm div}_{\wt \Om} 
\, \a = 0$.}

\medskip

\noindent {\it Proof.} It is easy to see that for an arbitrary 
function $f$ on $M$ one has:
\begin{equation}
{\rm div}_{f\Om} \, \g = {\rm div}_{\Om} \, \g + \frac{\g \, 
(df)}{f} \label{eq23}
\end{equation}
for any polyvector field $\g$. In the case when $\g = \a$, the 
Poisson bivector field, at each point $x \in M$ we have a map
$$
\ov{\a}_x : T_x^* \, M \ra T_x \, M
$$
defined by the formula
\begin{equation}
\ov{\a} \, (\om) = i_{\om} \, \a. \label{eq24}
\end{equation}

By definition, the image ${\rm Im} \, \ov{\a}_x$ is the tangent 
space to the symplectic leaf passing through a point $x$. The 
kernel ${\rm Ker} \, \ov{\a}_x$ consists of the 1-forms 
vanishing on the symplectic leaf.

When $f = \psi$ is constant along each symplectic leaf, the 
1-form $df$ vanishes on the tangent space to each symplectic 
leaf, and, therefore, $\a (d\psi) = 0$.~$\Box$

\bigskip

\subsection{}\label{ssec3.2} Now we consider the 
Poisson-Kirillov bivector field $\a$ on the space $\hbox{\goth 
g}^*$ dual to a semisimple Lie algebra. By 
Theorem~\ref{ssec2.1}, a constant volume form $\Om$ satisfies 
condition ${\rm div}_{\Om} \, \a = 0$. Then, by 
Key-lemma~\ref{ssec3.1}, we have
\begin{equation}
{\rm div}_{\psi \Om} \, \a = 0 \ \hbox{for any smooth $\psi$ 
constant on symplectic leaves.} \label{eq25}
\end{equation}
Finally, by Theorem~\ref{sssec1.3.1},
\begin{equation}
\int_{\hbox{\goth g}^*} f \cdot g \cdot \psi \cdot \Om = 
\int_{\hbox{\goth g}^*} (f * g) \cdot \psi \cdot \Om 
\label{eq26}
\end{equation}
for any $\psi$ constant along the symplectic leaves. In this 
way, we obtain a set of linear functionals on $C_{\rm 
comp}^{\ify} (\hbox{\goth g}^*)$ with the same values on $f 
\cdot g$ and $f * g$ (here $f,g \in C_{\rm comp}^{\ify} 
(\hbox{\goth g}^*))$. We want to deduce from this fact that $f 
\cdot g = f * g$ for $f,g \in [S^{\cdot} (\hbox{\goth 
g})]^{\hbox{\goth g}}$.

\subsubsection{Orbit method and Weyl unitary 
trick}\label{sssec3.2.1}

Since (\ref{eq26}) is stated for functions with compact support, 
it would be very useful to suppose that the symplectic leaves 
are compact. It is not true for semisimple algebras over $\Cb$. 
Fortunately, any such algebra has a {\it compact form}, i.e. 
such semisimple Lie algebra $\hbox{\goth g}^{\Rb}$ over $\Rb$ 
that $\hbox{\goth g}^{\Rb} \, \ot \, \Cb \simeq \hbox{\goth g}$ and 
that there exists a {\it compact} Lie group $G^{\Rb}$ such that 
Lie $G^{\Rb} = \hbox{\goth g}^{\Rb}$. By the Kirillov's theorem, 
the symplectic leaves are orbits of the coadjoint action of 
$G^{\Rb}$ on $(\hbox{\goth g}^{\Rb})^*$, and therefore, they are 
compact.

In the course of the proof of Theorem~\ref{ssec1.3} we 
will suppose that the symplectic leaves are compact.

\subsubsection{}\label{sssec3.2.2} When the symplectic leaves 
are compact, one can define a set of traces on the algebra 
$(C^{\ify} (\hbox{\goth g}^*) , *)$ (and $(S^{\cdot} (\hbox{\goth 
g}) , *)$), not only on the algebra $(C_{\rm comp}^{\ify} 
(\hbox{\goth g}^*) , *)$.

\medskip

\noindent {\bf Lemma.} {\it Let $\psi$ be a function on 
$\hbox{\goth g}^*$ with a} compact support {\it constant on each 
symplectic leaf, and $\Om$ be a volume form on $\hbox{\goth 
g}^*$ such that ${\rm div}_{\Om} \, \a = 0$ where $\a$ is the 
Poisson-Kirillov bivector field on $\hbox{\goth g}^*$. Then the 
functional
\begin{equation}
f \mpo \int_{\hbox{\goth g}^*} f \cdot \psi \cdot \Om 
\label{eq27}
\end{equation}
is a trace on the algebra $(C^{\ify} (\hbox{\goth g}^*), *)$, and 
on $(S^{\cdot} (\hbox{\goth g}) , *)$.}

\medskip

\noindent {\it Proof.} When $\psi$ has a compact support, the form 
$\psi \cdot \Om$ is {\it not} a volume form, and, therefore, one 
can not apply the theory of \cite{FSh} to it. But the function 
$\psi$ can be represented as the limit of a sequence of functions 
$\{ \psi_n \}_{n \geq 1}$, non-vanishing on $\hbox{\goth g}^*$. We 
know that the functional $f \mpo \int_{\hbox{\goth g}^*} f \cdot 
\psi_n \cdot \Om$ is a trace on $(C_{\rm comp}^{\ify} (\hbox{\goth 
g}^*) , *)$ for each $n$, and, therefore, the limit of the 
sequences of the functionals $f \mpo \int_{\hbox{\goth g}^*} f 
\cdot \psi \cdot \Om = \build \lim_{n \ra \ify}^{} \int_{\hbox{\goth 
g}^*} f \cdot \psi_n \cdot \Om$ is a trace on $(C^{\ify} 
(\hbox{\goth g}^*) , *)$.~$\Box$

\bigskip

\subsection{Theorem.}\label{ssec3.3} {\it Let {\goth g} be a 
semisimple compact Lie algebra over $\Rb$. Then for any two 
polynomials $f,g \in S^{\cdot} (\hbox{\goth g})$ one has:
\begin{equation}
f \cdot g = f * g + \om \label{eq28}
\end{equation}
where $\om \in [S^{\cdot} (\hbox{\goth g}) , S^{\cdot} 
(\hbox{\goth 
g})]_*$. (Here $[ \ , \ ]_*$ denotes the commutant of the deformed 
algebra.)}

\medskip

We prove this theorem in 
Section~\ref{sssec3.3.1}--\ref{sssec3.3.3}.

\subsubsection{Lemma.}\label{sssec3.3.1} {\it Let us suppose that a 
Lie 
algebra {\goth g} is such that the symplectic leaves in 
$\hbox{\goth g}^*$ are compact, and such that there exists a 
volume 
form $\Om$ on $\hbox{\goth g}^*$ satisfying the condition ${\rm 
div}_{\Om} \, \a = 0$, where $\a$ is the Poisson-Kirillov bivector 
field on $\hbox{\goth g}^*$. Then
$$
[C^{\ify} (\hbox{\goth g}^*)]^{\hbox{\goth g}} \cap [C^{\ify} 
(\hbox{\goth g}^*) , C^{\ify} (\hbox{\goth g}^*)]_* = 0
$$
($[ \ , \ ]_*$ denotes the commutant of the deformed algebra).}

\medskip

\noindent {\bf Corollary.} {\it In the assumptions of the Lemma one 
has:
$$
[S^{\cdot} (\hbox{\goth g})]^{\hbox{\goth g}} \cap [S^{\cdot} 
(\hbox{\goth g}) , S^{\cdot} (\hbox{\goth g})]_* = 0 \, .
$$}

\noindent {\it Proof of the Lemma.} Let $f \in [C^{\ify} 
(\hbox{\goth g}^*)]^{\hbox{\goth g}}$; then $f$ is constant on 
each symplectic leaf, because the leaves are orbits of the 
coadjoint 
action of a Lie group such that Lie $G = \hbox{\goth g}$ on 
$\hbox{\goth g}^*$. If $f$ is not zero function, one can choose a 
function $\psi$ on $\hbox{\goth g}^*$ with a compact support and 
such that
\begin{equation}
\int_{\hbox{\goth g}^*} f \cdot \psi \cdot \Om \ne 0 \, . 
\label{eq29}
\end{equation}
On the other hand, the functional
$$
f \mpo \int_{\hbox{\goth g}^*} f \cdot \psi \cdot \Om
$$
is a trace on $(C^{\ify} (\hbox{\goth g}^*) , *)$ according to 
Lemma~\ref{sssec3.2.2}. Therefore, if $f = [g_1 , g_2]_*$, then 
$\int_{\hbox{\goth g}^*} f \cdot \psi \cdot \Om = 0$, in a 
contradiction with (\ref{eq29}).~$\Box$

\bigskip

\subsubsection{Lemma.}\label{sssec3.3.2} {\it Let {\goth g} be a 
semisimple Lie algebra. Then the composition
$$
[S^{\cdot} (\hbox{\goth g})]^{\hbox{\goth g}} \hra S^{\cdot} 
(\hbox{\goth g}) \ra S^{\cdot} (\hbox{\goth g}) \bigl/ [S^{\cdot} 
(\hbox{\goth g}) , S^{\cdot} (\hbox{\goth g})]_*
$$
is an isomorphism.}

\medskip

\noindent {\it Proof.} According to Corollary~\ref{sssec3.3.1} and 
discussion in Section~\ref{sssec3.2.1} we know that
\begin{equation}
[S^{\cdot} (\hbox{\goth g})]^{\hbox{\goth g}} \cap [S^{\cdot} 
(\hbox{\goth g}) , S^{\cdot} (\hbox{\goth g})]_* = 0 \label{eq30}
\end{equation}
for any semisimple Lie algebra. As a {\goth g}-module (with respect 
to 
the 
adjoint action
in the deformed algebra), each filtration component $S^{\cdot} 
(\hbox{\goth g})_N$ of the space $S^{\cdot} (\hbox{\goth g})$ 
decomposes 
in the direct sum of finite-dimensional irreducible 
representations 
corresponding to dominant highest weights. We have:
\begin{equation}
S^{\cdot} (\hbox{\goth g})_N = \bigoplus_{{\rm dominant} \, \lb} 
L_{\lb}^{\op n_{\lb}} \, . \label{eq31}
\end{equation}
In these notations $[S^{\cdot} (\hbox{\goth g})_N]^{\hbox{\goth g}} = 
L_0^{\op n_0}$. On the other hand, consider the space $S^{\cdot} 
(\hbox{\goth g}) / [S^{\cdot} (\hbox{\goth g}) , S^{\cdot} 
(\hbox{\goth g})]_*$ as a {\goth g}-module.  Again, it is a trivial 
module,
$$
S^{\cdot} (\hbox{\goth g})_N \bigl/ ([S^{\cdot} (\hbox{\goth g}) , 
S^{\cdot} (\hbox{\goth g}) ]_*) \cap S^{\cdot} (\hbox{\goth g})_N
= L_0^{\op m_0} \, .
$$
It follows from (\ref{eq30}) and from the semisimplicity of 
$S^{\cdot} (\hbox{\goth g})$ as a {\goth g}-module that $n_0 = 
m_0$.~$\Box$

\bigskip

\subsubsection{}\label{sssec3.3.3} Here we conclude the proof of 
Theorem~\ref{ssec3.3}. It follows from Lemma 3.3.2 that
$f \cdot g - f * g = a + l$ where
$a \in  [S^{\cdot} (\hbox{\goth g}) , S^{\cdot} (\hbox{\goth g})]_*$
and $l \in [S^{\cdot} 
(\hbox{\goth g})]^{\hbox{\goth g}}.$
We have:
$$\int_{\hbox{\goth g}^*} (f \cdot g - f * g) \cdot \psi \cdot \Om = 
0$$
for any $\psi$, as well as 
$$
\int_{\hbox{\goth g}^*} a \cdot \psi \cdot \Om = 
0.$$
On the other hand, the function $l$ is constant along
the symplectic leaves, and one can choose a function $\psi$
with a compact support such that $\int_{\hbox{\goth g}^*}
 l \cdot \psi \cdot \Om \ne 0 \, .$
 Therefore, $l = 0$.

Theorem~\ref{ssec3.3} is proven.~$\Box$

\bigskip

\subsection{Theorem.}\label{ssec3.4} {\it For a semisimple Lie 
algebra {\goth g} and for $f,g \in [S^{\cdot} (\hbox{\goth 
g})]^{\hbox{\goth g}}$ one has:
$$
f \cdot g = f * g \, .
$$
}

\medskip

\noindent {\it Proof.} We know from Theorem~\ref{ssec3.3} that $f 
\cdot g = f * g + \om$. But when $f,g \in [S^{\cdot} (\hbox{\goth 
g})]^{\hbox{\goth g}}$, $f \cdot g$ and $f * g$ also belong to 
$[S^{\cdot} (\hbox{\goth g})]^{\hbox{\goth g}}$, and, therefore, 
$\om \in [S^{\cdot} (\hbox{\goth g})]^{\hbox{\goth g}}$. On the 
other hand, $\om \in [S^{\cdot} (\hbox{\goth g}) , S^{\cdot} 
(\hbox{\goth g})]_*$ by Theorem~\ref{ssec3.3}. Then $\om = 0$ by 
Corollary~\ref{sssec3.3.1}.

\subsection{Lemma~\ref{sssec3.3.1} and Harish-Chandra 
isomorphism}\label{ssec3.5}

Although the result of Lemma~\ref{sssec3.3.1} seems to be new, the 
Corollary~\ref{sssec3.3.1} is well-known. Its classical proof is 
based on the Harish-Chandra isomorphism (\cite{D}, Ch.~7.4). In 
this section we discuss this proof and relations between the 
Harish-Chandra isomorphism and Theorem~\ref{sssec1.3.1}.

We prove here the result of Corollary~\ref{sssec3.3.1} in its 
classical form: Let {\goth g} be a semisimple Lie algebra; then $Z 
\, (\Uc (\hbox{\goth g})) \cap [\Uc (\hbox{\goth g}) , \Uc 
(\hbox{\goth g})] = 0$. (Here $\Uc (\hbox{\goth g})$ is the 
universal enveloping algebra of the Lie algebra {\goth g} and $Z 
\, (\Uc (\hbox{\goth g}))$ is its center.)

The Harish-Chandra isomorphism is a map of algebras
$$
\vp_{\rm HC} : Z \, (\Uc (\hbox{\goth g})) \ra \Cb \, [\hbox{\goth 
h}]^W
$$
where {\goth h} is a Cartan subalgebra in {\goth g}, $W$ is the 
Weyl 
group. The algebra $\Cb \, [\hbox{\goth h}]^W$ consists of all 
polynomials $p : \hbox{\goth h}^* \ra \Cb$ such that
\begin{equation}
p \, (\lb) = p \, (w \, (\lb + \rho) - \rho) \label{eq32}
\end{equation}
for any $w \in W$ and any $\lb \in \hbox{\goth h}^*$; $\rho \in 
\hbox{\goth h}^*$ is the half-sum of the positive roots.

This map is constructed as follows: for any $\lb \in \hbox{\goth 
h}^*$ let $M_{\lb}$ be the Verma module with the highest weight 
$\lb$. Any element $z \in Z \, (\Uc (\hbox{\goth g}))$ acts on 
$M_{\lb}$ by multiplication on a constant. In such a way we 
obtain a map $\wt{\vp}_{\rm HC} : Z \, (\Uc (\hbox{\goth g})) \ra 
\Cb \, [\hbox{\goth h}]$. Now, if $M_{\lb'} \sbs M_{\lb}$ for some 
$\lb , \lb' \in \hbox{\goth h}^*$, the one should have 
$\wt{\vp}_{\rm HC} (z) (\lb) = \wt{\vp}_{\rm HC} (z) (\lb')$ for 
any $z \in Z \, (\Uc (\hbox{\goth g}))$. The classical result of 
I.N.~Bernstein, I.M.~Gelfand and S.I.~Gelfand (\cite{D}, Ch.~7.6) 
states that for a dominant $\lb \in \hbox{\goth h}^*$ the Verma 
module $M_{w (\lb + \rho) - \rho}$ is a submodule of $M_{\lb}$ for 
any $w \in W$. Therefore, $\wt{\vp}_{\rm HC} (z) \in \Cb \, 
[\hbox{\goth h}]^W$ for any $z \in Z \, (\Uc (\hbox{\goth g}))$, 
and we obtain a map
$$
\vp_{\rm HC} : Z \, (\Uc (\hbox{\goth g})) \ra \Cb \, [\hbox{\goth 
h}]^W \, .
$$
The Harish-Chandra theorem states that it is an isomorphism.

Now we are going to deduce from this theorem that $Z \, (\Uc 
(\hbox{\goth g})) \cap [\Uc (\hbox{\goth g}) , \Uc (\hbox{\goth 
g})] = 0$. Indeed, for a dominant weight $\lb \in \hbox{\goth 
h}^*$ 
there exists the unique irreducible {\it finite-dimensional} 
{\goth 
g}-module $L_{\lb}$ with highest weight $\lb$. The module 
$L_{\lb}$ 
is a quotient of the Verma module $M_{\lb}$, and, therefore, any 
element $z \in Z \, (\Uc (\hbox{\goth g}))$ acts on $L_{\lb}$ in 
the same way as it does on $M_{\lb}$. Hence, for a dominant $\lb$,
$$
\vp_{\rm HC} (z) (\lb) = \frac{1}{\dim \ L_{\lb}} \ {\rm 
Tr}_{L_{\lb}} (z) \, .
$$
Now it is clear that $Z \, (\Uc (\hbox{\goth g})) \cap [\Uc 
(\hbox{\goth g}) , \Uc (\hbox{\goth g})] = 0$, because for any $z$ 
from the intersection one has $\vp_{\rm HC} (z) (\lb) = 0$ for 
dominant $\lb$, and, therefore, for all $\lb \in \hbox{\goth 
h}^*$. 
Then $z = 0$ because $\vp_{\rm HC}$ is an isomorphism.

Philosophically, a version of the Harish-Chandra map may exist for 
any Lie algebra {\goth g}, not necessarily semisimple, or, more 
generally, for any Poisson structure on a vector space $V$. It is 
not true literally, but the well-known application of this 
philosophy 
is the Kirillov's work on Duflo formula. In this work the Duflo 
isomorphism was obtained using the Harish-Chandra map, and was 
rewritten in a form independent of a choice of a Cartan subalgebra. 
Then the obtained formula make sense for any 
finite-dimensional Lie algebra, not only semisimple.

>From this point of view, the statement of Lemma~\ref{sssec3.3.1} is 
another incarnation of this idea. Our proof, based on the result 
of 
\cite{FSh}, does not use any special properties of semisimple Lie 
algebras except the compactness of the orbits. It can be easily 
generalized as follows: consider a vector spase $V$
(or, more generally, an oriented manifold $M$)
with a Poisson structure $\a$ with compact symplectic leaves, suppose 
also that there exists a volume form $\Om$ on $V$ such that ${\rm 
div}_{\Om} \, \a = 0$. Denote by $[C^{\ify} (V)]^{\a}$ the 
functions which are constant along the symplectic leaves. Then
$$
[C^{\ify} (V)]^{\a} \cap [C^{\ify} (V) , C^{\ify} (V)]_* = 0 \, .
$$

\subsection{Vanishing of the wheels}\label{ssec3.6}

Let $k$ be a minimal integral number such that $W_{2k} \ne 0$. 
Then, by Theorem~\ref{ssec1.2}~(i) and by Theorem~\ref{ssec3.4}, 
the element ${\rm Tr}_{2k} \in [S^{2k} (\hbox{\goth 
g}^*)]^{\hbox{\goth g}}$ should be a derivation of the algebra 
$[S^{\cdot} (\hbox{\goth g})]^{\hbox{\goth g}}$ for any semisimple 
Lie algebra {\goth g}. One can show that it is not true for Lie 
algebra $s\ell_n$ for sufficiently large $n$. (The same argument was 
used by Maxim Kontsevich in his proof of Duflo theorem, see 
\cite{K}, Sect.~8.3.4.) Theorem~\ref{ssec1.3} is proven.

\subsection{Remarks.}\label{ssec3.7} 1. The integrals, corresponding
to the wheels, do not 
depend on a choice of the angle function (\cite{K}, Sect.~6.2 and 
Sect.~8) We 
have proved that $W_{2k} = 0$ for the harmonic angle function (since 
the
 result of \cite{FSh} holds for harmonic angle function).
 Hence, it is true for any 
angle function.

\smallskip

\noindent 2. The odd wheels $W_{2k+1}$ are equal to zero by a much 
more simple reasons: it follows from the fact that for the 
harmonic angle function the propagator is invariant under the map 
$p : z \mpo - \ov z$ ($p : \Hc \ra \Hc$, where $\Hc$ is the 
complex upper half-plane), as it was noted in \cite{K}. The same 
argument shows that $B_k (f,g) = (-1)^k \, B_k (g,f)$ for $k \geq 
0$, {\it any} $f,g \in C^{\ify} (V)$ and {\it any} Poisson 
bivector field (not necessarily linear) on a vector space $V$ (see 
\cite{CF}).